\documentclass{article}
\usepackage[utf8]{inputenc}
\usepackage{xcolor}
\usepackage{amsfonts}
\usepackage{amsmath}
\usepackage{appendix}
\usepackage{graphicx}
\usepackage{subfig}
\numberwithin{equation}{section}

\title{A caustic terminating at an inflection point }
\author{J. R. Ockendon, H. Ockendon\\Mathematical Institute, University of Oxford\\R. H. Tew\\School of Mathematical Sciences, University of Nottingham\\
D. P. Hewett, A. Gibbs\\Department of Mathematics, University College, London}
\begin{document}
\maketitle

\begin{abstract}
We present an asymptotic and numerical study of the evolution of an incoming wavefield which has a caustic close to a curve with an inflection point.  Our results reveal the emergence of a wavefield which resembles that of a shadow boundary but has a maximum amplitude along the tangent at the inflection point.
\end{abstract}
\section{Introduction}
In this brief paper we present an asymptotic and numerical study of a two-dimensional wavefield with a caustic near the cubic curve $y+\frac{1}{3}\gamma x^3=0$ as $x\rightarrow -\infty$ in the far field. 
%an incoming caustic that is localised in the $(x,y)$-plane near 
%a caustic that is close to the curve ${y+\frac{1}{3}\gamma x^3=0}$ as $x\rightarrow -\infty$ in the far field. 
The problem was introduced in \cite{ocktew21} along with an explicit integral representation in \cite[(4.7)]{ocktew21} for the wavefield in a small neighbourhood of the origin, where the cubic has an inflection point and the caustic terminates. Unfortunately, \cite[(4.7)]{ocktew21} contained some errors, but also gave little insight into the outgoing field as $x\rightarrow +\infty$. This paper both corrects the errors in \cite[(4.7)]{ocktew21} and provides a detailed asymptotic description and a numerical solution of the outgoing wavefield.

The main rationale for this paper is that improved understanding of this wavefield may give intuition about the far-field solution of the famous Popov inflection point problem, 
which 
%which, in the context of wave scattering in domains with boundaries, 
describes the scattering of an incident whispering gallery wavefield propagating along a concave portion of a scatterer boundary when it reaches an inflection point of the boundary. This canonical problem, which remains unsolved in closed form, has a lengthy history, the most recent review of which is \cite{smykam22}. 
%In its simplest form, the 
The Popov problem seeks a wavefield in the region $y>-\frac{1}{3}\gamma x^3$ satisfying either a Dirichlet or Neumann boundary condition on the curve $y+\frac{1}{3}\gamma x^3=0$, and tending to an Airy function close to this curve as $x\to-\infty$. In this paper we consider the evolution of an Airy function wavefield past the inflection point of the cubic in the \textit{absence} of any boundary. 

In \S\ref{asymptotic} we introduce appropriate local curvilinear coordinates $\left(S,N\right)$ based on the cubic curve, and then derive a representation for the wavefield close to the inflection point in the form of an integral involving a complicated phase function of $S$ and $N$, correcting \cite[(4.7)]{ocktew21}. We then consider the outgoing field in the limit as $S\rightarrow \infty$, which we show to have a transition between bright and dark in the vicinity of the positive $x$-axis.
In \S\ref{numerics} we reformulate the solution in terms of integrals studied in \cite{hewocksmy19}, which are appropriate for the numerical scheme introduced in \cite{gibhewhuy23}, allowing us to present a validation of the asymptotic analysis.
Finally, in \S\ref{conclusions} we offer some conclusions. 

\section{Asymptotic behaviour of the wavefield near the inflection point}\label{asymptotic}
We align a two-dimensional Cartesian frame so that, for sufficiently small values of $x$ and $y$, the caustic has equation $y+\frac{1}{3}\gamma x^3=0$, $x<0$, for some positive constant $\gamma$. Local scalings on these coordinates, which both preserve this cubic profile exactly and allow us to model wave propagation through the inflection point at the origin are found to be $x=k^{-1/5}X$, $y=k^{-3/5}Y$, with $X,Y=O(1)$. One can then seek an approximate solution to the Helmholtz equation $(\Delta + k^2)\phi = 0$ of the form $\phi \sim e^{ikx}\tilde {A}\left(X,Y\right)$ as $k\to\infty$, with ${\tilde{A}}$ satisfying the parabolic wave equation 
%(see e.g.\ \cite{hewocksmy19})
\begin{equation} 
\label{eq: pwe}
\frac{\partial ^2\tilde{A}}{\partial Y^2}+2i\frac{\partial \tilde{A}}{\partial X}=0.
\end{equation}

An alternative approach, adopted and reviewed in \cite{ocktew21}, is to switch to local curvilinear coordinates intrinsic to the curve ${y+\frac{1}{3}\gamma x^3=0}$ involving arc-length $s$ and normal distance $n$ along and from it, respectively. It then follows that, correct up to and including terms of ${O\left(k^{-1}\right) }$ as $k\rightarrow \infty$,
\[x\sim s+\gamma ns^2-\frac{1}{10}\gamma^2s^5,\quad y\sim n-\frac{1}{3}\gamma s^3.\]
Thus, with the scalings ${s=k^{-1/5}S,\,n=k^{-3/5}N }$ $\left(|S|,\,|N|=O(1)\right)$ and seeking an alternative local solution ${\phi \sim e^{iks}A\left(S,N\right) }$, we have that, to lowest order,
\begin{align*}
ikx&=iks+i\left(\gamma NS^2-\frac{1}{10}\gamma ^2S^5\right), \\
\tilde{A}\left(X,Y\right)&=e^{-i\left(\gamma NS^2-\frac{1}{10}\gamma ^2S^5\right)}\,A\left(S,N\right),\\
X&=S,\quad Y=N-\frac{1}{3}\gamma S^3. 
\end{align*}
%\begin{itemize}
%    \item[(a)] ${ikx=iks+i\left(\gamma NS^2-\frac{1}{10}\gamma ^2S^5\right) }$; 
%    \item[(b)] ${\tilde{A}\left(X,Y\right)=e^{-i\left(\gamma NS^2-\frac{1}{10}\gamma ^2S^5\right)}\,A\left(S,N\right) }$; 
%    \item[(c)] ${X=S,\,Y=N-\frac{1}{3}\gamma S^3 }$. 
%\end{itemize}    
It then follows that (\ref{eq: pwe}) can be re-written in terms of $S$ and $N$ as
\[\left(\frac{\partial ^2}{\partial N^2}+2i\left(\frac{\partial}{\partial S}+\gamma S^2\frac{\partial}{\partial N}\right)\right)\left(e^{-i\left(\gamma NS^2-\frac{1}{10}\gamma ^2S^5\right)}A\left(S,N\right)\right)=0,\]
which leads to the Popov equation
\begin{equation}
    \label{eq: popoveqn}
    \frac{\partial ^2A}{\partial N^2}+2i\frac{\partial A}{\partial S}+4\gamma NSA=0.
\end{equation}
%as presented in \cite{pop79} and \textbf{(REF TL2)}, follows.\\
%As described in \textbf{REF: 2021}, when we work in curvilinear coordinates $(s,n)$, where $s$ is arclength and $n$ the normal distance from ${ y+\frac{1}{3}\gamma x^3=0}$, the incoming wavefield away from the origin takes the form ${e^{iks}A(s,n) }$ as $k\rightarrow \infty$, where $A$ is an Airy function. \\
%However, near the origin, we write
%\[ s=k^{-1/5}S,\,n=k^{-3/5}N\]
%which means that the lowest order Helmholtz equation becomes
%\begin{equation} 
%\label{eq: popoveqn}
%\frac{\partial ^2A}{\partial N^2}+2i\frac{\partial A}{\partial S}+4\gamma NSA=0.
%\end{equation}.

In order to match with the incoming Airy function we require that
\begin{equation}
\label{eq: Aincom}
    A\sim 2\pi \left(-4\gamma S\right)^{1/3}\mbox{Ai}\left[ \left(-4\gamma S\right)^{1/3}N\right], \qquad S\rightarrow -\infty,
\end{equation}
where the first ${ \left(-S\right)^{1/3}}$ term reflects the dependence of the incoming field on the curvature of the associated caustic as $S\rightarrow -\infty$ and the constants in the pre-factor multiplying the Airy function are inserted for later convenience.
Since the data is analytic, we expect the solution of (\ref{eq: popoveqn}) to be analytic for all $S,N$.

The Fourier transform procedure outlined in \cite{ocktew21} may now be used to show that the desired solution of (\ref{eq: popoveqn}) is
\begin{equation}
    \label{eq: solforpopov}
    A=\int_{-\infty}^{\infty}\,e^{-\frac{i}{2}E}\,d\lambda ,
\end{equation}
where
\begin{align}
\nonumber 
E&=S\left(\lambda+\gamma S^2\right)^2-\frac{2\gamma}{3}\left(\lambda +\gamma S^2\right)S^3+\frac{1}{5}\gamma ^2S^5\\
\nonumber
&\qquad +\frac{8}{15\gamma ^{1/2}}\left(\lambda+\gamma S^2\right)^{5/2}+2N\lambda ;
\end{align}
here, ${ \left(\lambda +\gamma S^2\right)^{1/2}}$ is positive when $\lambda +\gamma S^2>0$ and equal to ${ -i|\lambda +\gamma S^2|^{1/2}}$ when $\lambda +\gamma S^2<0$, and the integral is taken just above the negative real $\lambda$ axis and along the positive $\lambda$ axis. This equation corrects \cite[(4.7)]{ocktew21}.

In order to deal with the branch point at $\lambda = - \gamma S^2$, it is convenient to write
\begin{equation}\label{eq:Aleft}
A=\left(\int_{-\gamma S^2}^{\infty}\,+\int_{-\infty}^{-\gamma S^2}\right)\,e^{-\frac{i}{2}E}\, d\lambda =I_++I_-\,\,\mbox.
\end{equation}
This decomposition can be justified by deforming the integration path from one along $\Im{\lambda}>0$ and showing that the contribution from the region near $\lambda = -\gamma S^2$ is negligible as $\Im{\lambda} \downarrow 0$.  The remainder of this section will be based on (\ref{eq:Aleft}).

\subsection{The limit $S \rightarrow - \infty$}

Noting that when $S$ is large and negative, ${\left(\lambda +\gamma S^2\right)^{5/2}=-\gamma ^{5/2}S^5\left(1+\frac{\lambda}{\gamma S^2} \right)^{5/2}}$, we can expand $E$ for large negative $S$ in $I_+$. This gives that
${ E\sim 2N\lambda - \frac{\lambda^3}{6\gamma S} + \ldots}$ and a formal asymptotic expansion shows that the leading term in $I_+$ is
\begin{equation}
\nonumber
   \int_{-\gamma S^2}^{\infty}e^{-i(N\lambda - \frac{\lambda^3}{12\gamma S})}d\lambda, 
\end{equation}
which, by putting $\lambda = \left(-4\gamma S\right)^{1/3} \tau$, is asymptotic to 
\begin{eqnarray}
\nonumber
      (-4\gamma S)^{1/3}\int_{-\infty}^{\infty}e^{-i\left(\frac{{\tau}^3}{3}+(-4\gamma S)^{1/3}N\tau \right)}d\tau     
    = 2\pi \left(-4\gamma S\right)^{1/3}\mbox{Ai}\left( \left(-4\gamma S\right)^{1/3}N\right)
\end{eqnarray}
as $S\rightarrow -\infty$. Thus, since $I_-$ tends to zero as $S\rightarrow -\infty$, $A$ matches with (\ref{eq: Aincom}) as expected.  An estimate of the rate at which the Airy function is attained for large values of $-S$ would involve a complicated analysis which we will not attempt in this paper.

\subsection{The limit $S\rightarrow +\infty$}

The situation is more interesting when $S\rightarrow +\infty$, when it is convenient to write
$\lambda+\gamma S^2= -S^2 T$ for $T\ge 0$ in $I_-$
so that
\begin{equation}\label{eq:Iminus}
I_- = S^2e^{i\Phi}\int_0^{\infty}\, e^{-S^5\left(\frac{i}{2}\left(T^2-2KT\right)
+\frac{4}{15\gamma^{1/2}}T^{5/2}\right)}\,dT.     
\end{equation}
where ${\Phi=\gamma NS^2-\frac{\gamma^2 S^5}{10}}$ and 
\[K = \frac{N}{S^3}-\frac{\gamma}{3} \] 
is a measure of transverse distance from the $x$-axis.
Similarly, writing $\lambda+\gamma S^2= S^2 T$ for $T\ge 0$ in $I_+$ leads to
\begin{equation}\label{eq:Iplus}
I_+ = S^2e^{i\Phi}\int_0^{\infty}\, e^{-\frac{i}{2}S^5\left(T^2+2KT+\frac{8}{15{\gamma}^{1/2}}T^{5/2}\right)}\,dT.     
\end{equation}

\subsubsection{$K\ge0$, $K$ of $O(S^{-2})$ or larger}

We begin by considering $I_-$ which we see from (\ref{eq:Iminus}) has a stationary phase point at $T=K$. The contribution from this point to the value of the integral is exponentially small as long as $S^2K\gg 1$ and, in this case, the main contribution to $I_-$ comes from a region near $T=0$ where $T^{5/2}$ is negligible to lowest order. However the $T^2$ term in the exponent needs to be retained to ensure convergence of the integral.  Hence the leading term in the integral is
\begin{equation}
\nonumber
    \int_0^{\infty}\,e^{-\frac{i}{2}S^5\left(T^2-2KT\right)}\,dT = 
    \frac{1}{S^5K}\int_0^{\infty}e^{i\left(\bar{T} - \frac{\epsilon}{2}{\bar{T}}^2\right)}\,d\bar{T},
\end{equation}
where $\bar{T}= S^5KT$ and $\epsilon = \frac{1}{S^5K^2}$.
We now note that 
\begin{eqnarray}
\nonumber
    \int_0^{\infty}e^{i\left(\bar{T} - \frac{\epsilon}{2}{\bar{T}}^2\right)}\,d\bar{T}
     = \lim_{\delta \downarrow 0} \int_0^{\infty}\,e^{(i-\delta )\bar{T}}\left(1-\frac{i\epsilon {\bar{T}}^2}{2} + ...\right) d\bar{T}\\
\nonumber
     = \lim_{\delta \downarrow 0} \left( \frac{1}{\delta - i} - \frac{i\epsilon}{(\delta - i)^3} + ...\right)
     =i - \epsilon + ...   .
\end{eqnarray}     
The device of introducing the parameter $\delta$ is one way of deriving higher order terms in asymptotic expansions of stationary phase integrals, as described in \cite{olver74}.   It gives us the result that 
\begin{equation}\label{eq: IminusKplus}
    I_- \sim  e^{i\Phi}\frac{i}{S^3K}\left(1 + \frac{i}{S^5K^2} + ... \right),\quad S\rightarrow +\infty,\quad S^2K \gg 1,
\end{equation}
and it can also be used to obtain higher order corrections for the estimates we will derive in the rest of this paper.

The result corresponding to (\ref{eq: IminusKplus}) for $I_+$ can be read off by replacing $K$ by 
$-K$ and it tells us that $ A = I_+ + I_- = o(S^{-8}K^{-3})$ when $S\rightarrow + \infty$ with $K\geq 0$ and $S^2K\gg 1$.

None of these results applies when $K$ is of $O(S^{-2})$ or smaller.  In this regime, we write
$K= S^{-2}\hat{K}$ and $T=S^{-2}t$, so that 
\begin{equation}\label{eq:Iminus2}
    I_- = e^{i\Phi}\int_0^{\infty}\, e^{-S\frac{i}{2}\left(t^2-2\hat{K}t\right)
-\frac{4}{15\gamma^{1/2}}t^{5/2}}\,dt.  
\end{equation}
The dominant contribution as $S\rightarrow +\infty$ now comes from the stationary-phase point at 
$t=\hat{K}$ and we soon find $I_-$ is given by
\begin{equation}\label{eq:Iminus3}
    I_- \sim e^{i\Phi}\left(\frac{2\pi}{S}\right)^{1/2}e^{\left(\frac{iS{\hat{K}}^2}{2}-\frac{4{\hat{K}}^{5/2}}{15{\gamma}^{1/2}} - \frac{i\pi}{4}\right)}.
\end{equation}
Moreover, the use of the regularisation leading to (\ref{eq: IminusKplus}) reveals that the next order term in the stationary phase expansion of $I_-$ is of relative order $S^{-1}$.
However, there is also a contribution from the end point $t=0$, which can be estimated as in (\ref{eq: IminusKplus}) to give a leading order term 
${\frac{e^{-i\Phi}}{iS\hat{K}}}$.  This can be shown to cancel with a corresponding term in $I_+$ in which there is no stationary phase contribution.  
Hence we conclude that, with an expected relative error of $O(S^{-1})$,
\begin{equation}\label{eq:modAKplus}
  |A| \sim  \left(\frac{2\pi}{S}\right)^{1/2}e^{\left(-\frac{4{\hat{K}}^{5/2}}{15{\gamma}^{1/2}}\right)} 
\end{equation}
as $S\rightarrow +\infty$ with $K\ge 0$ and $\hat{K} = O(1)$.

We note that the derivation of (\ref{eq:Iminus3}) from (\ref{eq:Iminus2}) involves the assumption that the stationary phase point $\hat{K}$ is sufficiently far from the origin  that $\hat{K}\gg O(S^{-1/2})$. We will examine the region where $\hat{K} = O(S^{-1/2})$ shortly.\\
 
\subsubsection{$K\leq0$, $|K|$ of $O(S^{-2})$ or larger}

The calculations above can be repeated to show that the contributions to both $I_+$ and $I_-$ near $T=0$ cancel at least to leading order as $S\rightarrow +\infty$.   Hence we need only consider $I_+$, which has a stationary phase contribution.

It is now more convenient to write $ \lambda + \gamma S^2 = S^2 {\tau}^2$ so that 
\begin{equation}\label{eq:Iplus3}
I_+ = 2S^2e^{i\Phi}\int_0^{\infty}\, e^{-\frac{i}{2}S^5g(\tau)}\, \tau\,d\tau, 
\end{equation} 
where 
\begin{eqnarray}
  g(\tau)=2K{\tau}^2+{\tau}^4+\frac{8}{15\gamma^{1/2}}{\tau}^5.
\end{eqnarray}
Hence there is always one stationary phase point in $\tau > 0$, $\tau = {\tau}_0(K)$ say, which is the positive root of 
\begin{equation}
 K + {\tau_0}^2 + \frac{2}{3{\gamma}^{1/2}}{\tau_0}^3 = 0  
\end{equation}
and is such that ${\tau}_0(0) = 0$ and ${\tau}_0 = O({|K|}^{1/3})$ as $K\rightarrow -\infty$.
The leading order term in the expansion of the integral in (\ref{eq:Iplus3}) is thus
\begin{equation}\label{eq;K<0I}
    e^{-\frac{i\pi}{4}}{\tau}_0e^{-\frac{i}{2}S^5g({\tau}_0)}
    \left(\frac{4\pi}{S^5 g''({\tau}_0)} \right)^{1/2},
\end{equation}
where ${g''\left( {\tau}_0\right)=8{{\tau}_0}^2\left(1+\frac{{\tau}_0}{\gamma ^{1/2}}\right) }$.   
Hence 
\begin{equation}\label{eq:modAKminus}
 |A| \sim \sqrt{\frac{2\pi}{S(1 + \frac{{\tau}_0}{\gamma^{1/2}})}}, \quad S\rightarrow +\infty,  \, K<0,
\end{equation}
and so, as $K\uparrow 0$, ${|A|\sim \frac{\sqrt{2\pi}}{S^{1/2}}}$, which is equal to (\ref{eq:modAKplus}) when $\hat{K}=0$.  However, the comments made after (\ref{eq:modAKplus}) apply equally to the derivation of (\ref{eq;K<0I}) which only holds if $S^{5/2}{\tau}_0 \gg 1$, and we will describe the inner layer when $K=O(S^{-5/2})$ in the next subsection. 

Meanwhile we note that the results (\ref{eq:modAKplus}) and (\ref{eq:modAKminus}) indicate that the far-field as $S\rightarrow +\infty$ is dominated by the region in which $K$ is negative and of O(1).
%, and this region will be referred to as the searchlight region.   
Plots of the numerical evaluation of $A$ will be given in the next section and comparison will be made with (\ref{eq:modAKplus}) and (\ref{eq:modAKminus}). 

\subsubsection{$K = O(S^{-5/2})$}

In this region, the sign of $K$ is no longer important and we will simply consider (\ref{eq:Iplus}) and (\ref{eq:Iminus}) when ${K=S^{-5/2}\bar{K}}$.  Writing ${T=S^{-5/2}\bar{t}}$ then, to leading order as $S\rightarrow +\infty$, 
\begin{eqnarray}
\nonumber 
    I_+&\sim & e^{i\Phi}S^{-1/2}\int_0^{\infty}\,e^{-i\bar{t}^2/2-i\bar{K}\bar{t}}\,d\bar{t}\\
    \label{eq:fresnel}
    &\sim& e^{i\Phi}S^{-1/2}e^{i\bar{K}^2/2}\int_{\bar{K}}^{\infty}\, e^{-iv^2/2}\,dv
\end{eqnarray}
where $v=\bar{K}+{\bar{t}}$.
Repeating the exercise for $I_-$ gives 
\begin{equation}
    I_-\sim e^{i\Phi}S^{-1/2}e^{i\bar{K}^2/2}\int_{-\bar{K}}^{\infty}\, e^{-iv^2/2}\,dv
\end{equation}
and hence 
\begin{equation}
\label{eq:constant}
I_+ + I_- \sim S^{-1/2}e^{i\Phi + i\bar{K}^2/2}
\int_{-\infty}^{\infty}\, e^{-iv^2/2}\,dv 
= \sqrt{2\pi}S^{-1/2}e^{i\Phi + i\bar{K}^2/2 -i\pi/4}.
\end{equation}  
Thus the amplitude $|A|$ in this region is independent of $\Bar{K}$ and matches with (\ref{eq:modAKplus}) as $\Bar{K}\rightarrow +\infty$ and with (\ref{eq:modAKminus}) as $\Bar{K} \rightarrow -\infty$.\\

Even though the amplitude of the far-field solution can be described analytically, we have only worked to the lowest order when obtaining the asymptotic expansions for $A$ and this has resulted in $A$ having discontinuous slope as a function of $K$.  
We expect the solution of (\ref{eq: popoveqn}) to be analytic everywhere and that  the asymptotic approximation will become increasingly smooth when taken to higher orders but, as mentioned after (\ref{eq: IminusKplus}), this is a challenging task.  
Hence it is helpful to compare these predictions with numerical calculations, which we do in the next section.

\section{Numerical validation} \label{numerics}
In this section we compare the asymptotic approximations obtained in the previous section with accurate numerical evaluations of the integral \eqref{eq: solforpopov}, obtained using the PathFinder software \cite{pathfinder}. This implements the algorithm described in \cite{gibhewhuy23}, which automates the numerical steepest descent method for oscillatory integrals (see e.g.\ \cite[\S5]{deahuyise18}), automatically performing appropriate contour deformations and dealing robustly with multiple coalescing stationary points. 

\begin{figure}[t!]
    \centering
    \includegraphics[width=.6\linewidth]{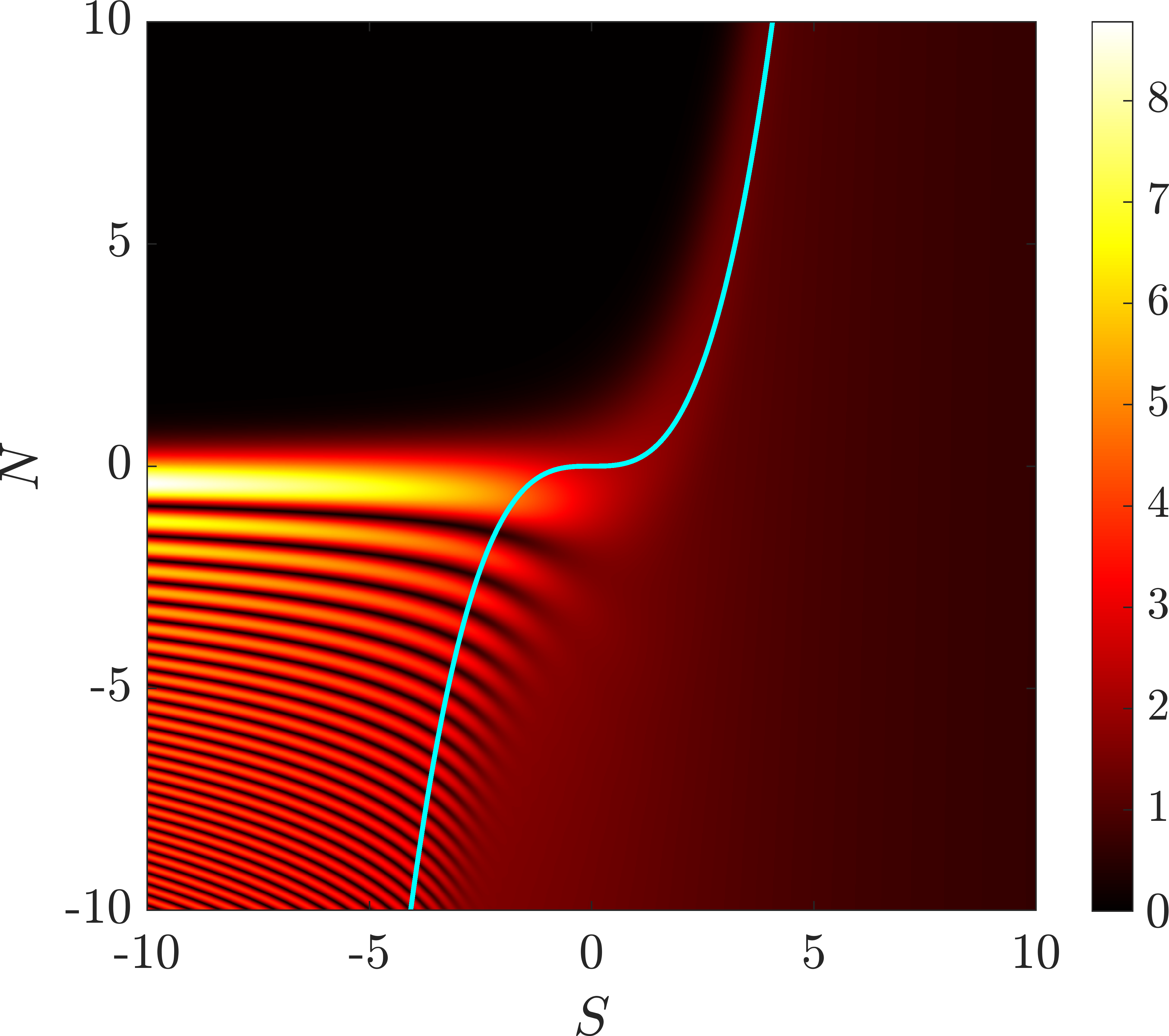}
    \caption{Plot of $|A|$ as a function of $S$ and $N$, along with the curve $K=0$, i.e.\ the cubic $N-\frac{\gamma}{3}S^3=0$, for $\gamma=4/9$.}
    \label{fig:1}
\end{figure}

\begin{figure}[t!]
    \centering
    \includegraphics[width=\linewidth]{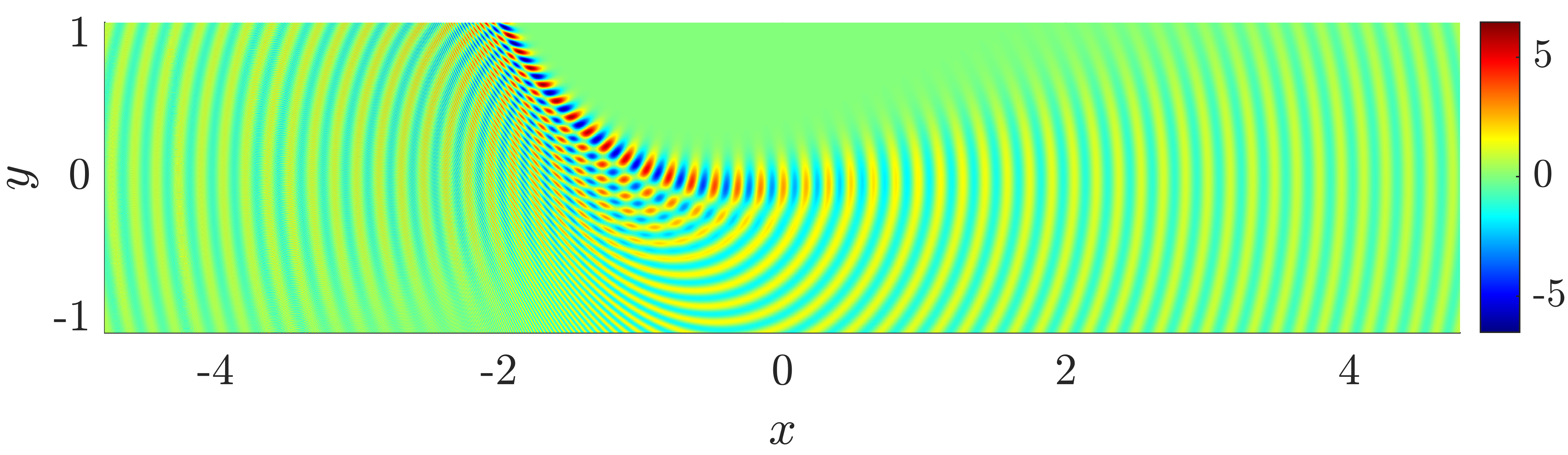}
    \caption{The real part of the approximate Helmholtz equation solution $-2e^{ikx}\tilde{A}_{32}(k^{1/5}x,k^{3/5}y)$ for $\gamma=4/9$ and $k=40$.}
    \label{fig:2}
\end{figure}

\begin{figure}[t!]
    \centering
    \subfloat[]{\includegraphics[height=50mm]{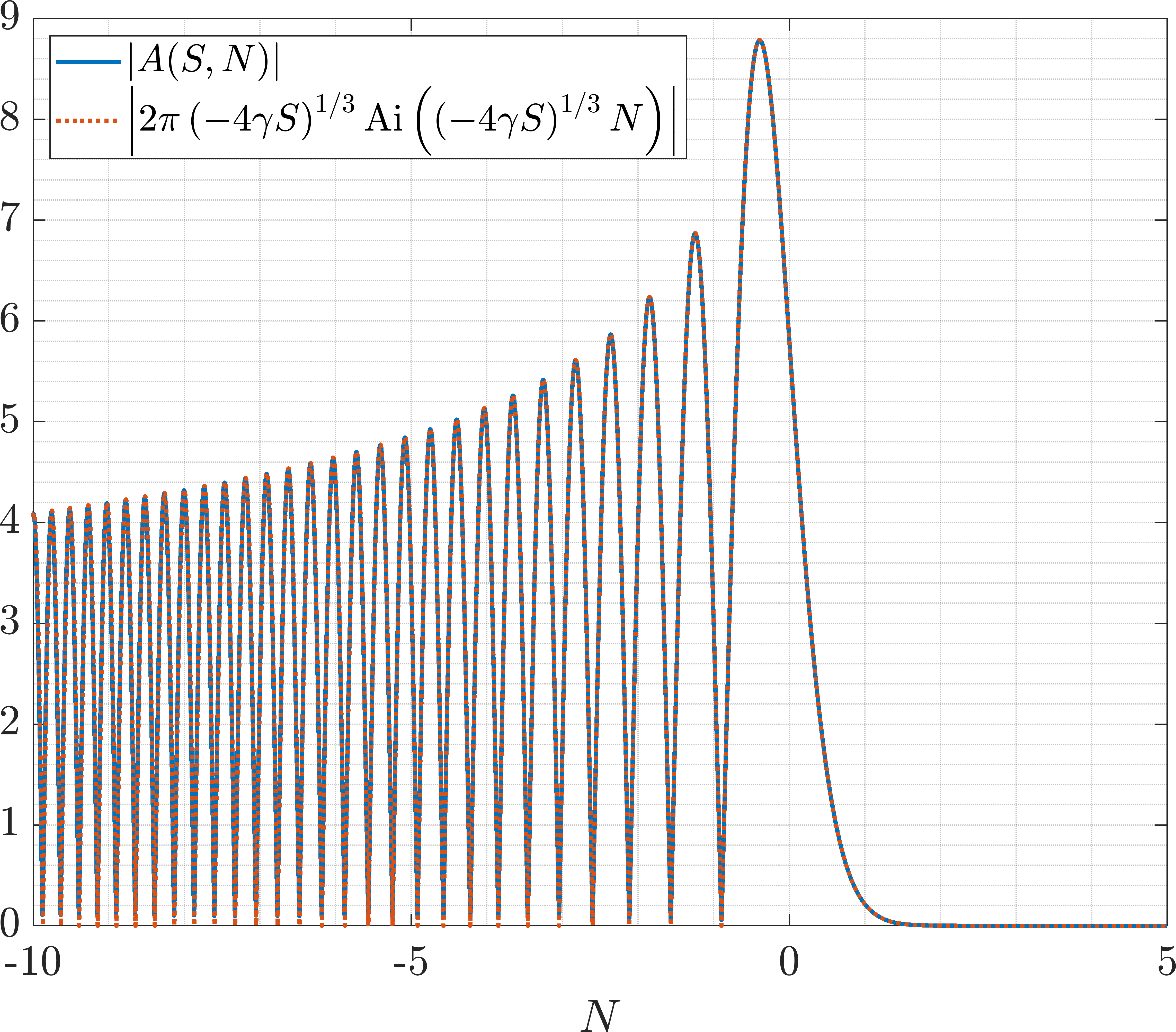}}
    \hspace{5mm}
    \subfloat[]{\includegraphics[height=50mm]{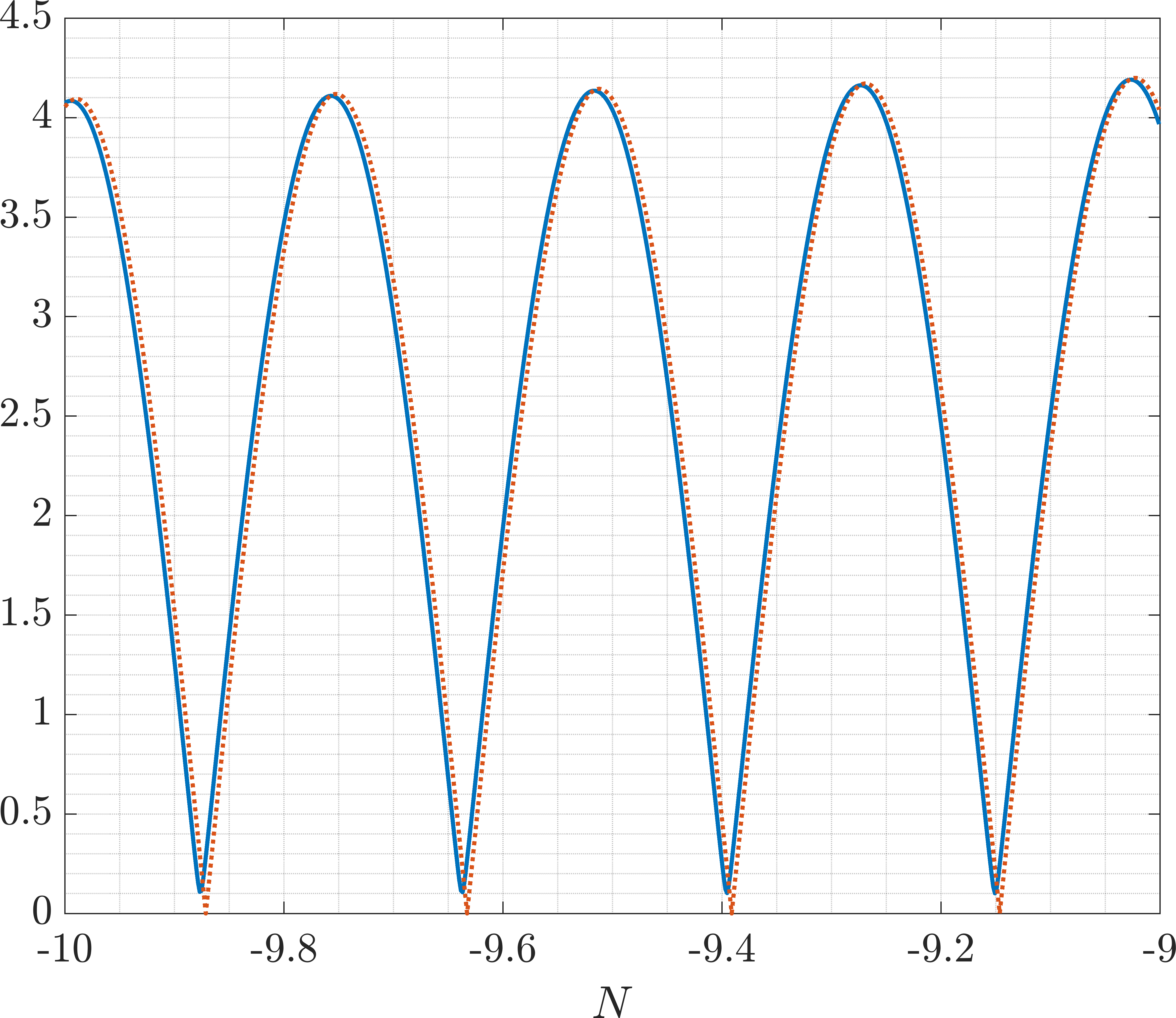}}
    \caption{(a) Plot of $|A|$ as a function of $N$ for fixed $S=-10$ and $\gamma=4/9$, showing the agreement with \eqref{eq: Aincom}. On this scale, the two curves are almost indistinguishable. (b) Zoom of (a) for $N\in[-10,-9]$.}
    \label{fig:3}
\end{figure}

\begin{figure}[t!]
    \centering
    \subfloat[$S=5$]{\includegraphics[height=49mm]{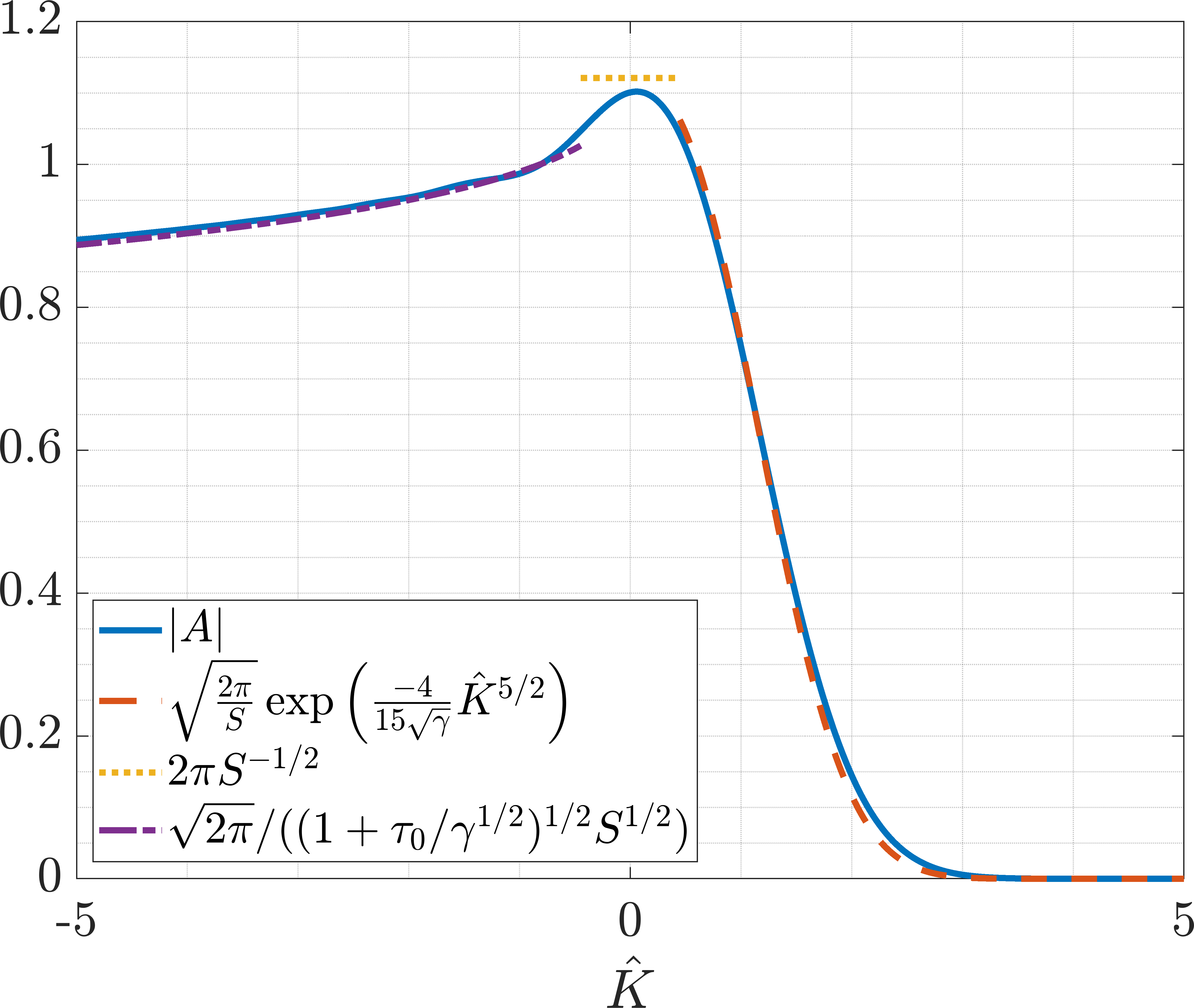}}
    \hspace{3mm}
    \subfloat[$S=10$]{\includegraphics[height=49mm]{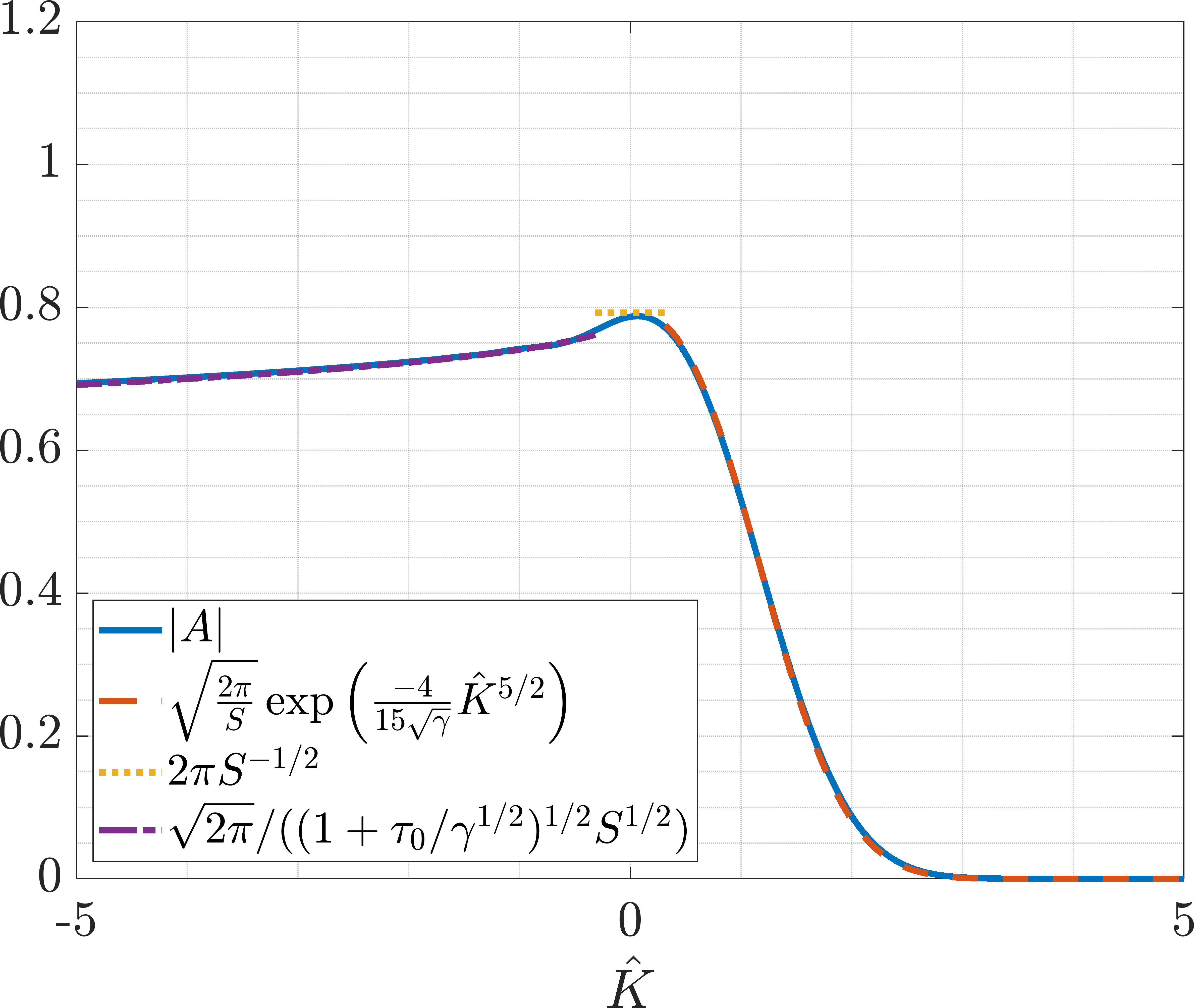}}
    %\subfloat[$S=10$]{\includegraphics[width=.5\linewidth]{A_cs_S_10_100qpts_10000ppts_zoom.png}}
    \caption{Plots of $|A|$ as a function of $\hat{K} = S^2K = \frac{N}{S}-\frac{\gamma}{3}S^2$ for $S=5$ and $S=10$, and $\gamma=4/9$, showing the comparison with \eqref{eq:modAKplus} (for $\hat{K}> 0$), \eqref{eq:modAKminus} (for $\hat{K}< 0$) and \eqref{eq:constant} (for $\hat{K}\approx 0$). 
    %(b) Zoom of (a) for $\hat{K}\in[0,5]$.}
    }
    \label{fig:4}
\end{figure}

To obtain an integral in a form amenable to evaluation using PathFinder, which in its current form applies only to integrals with polynomial phase, 
%which in its current form requires polynomial phase, 
we apply a change of variable to rewrite \eqref{eq: solforpopov} as
\begin{align}
\label{eq:A32}
A(S,N) =
%-\frac{1}{\pi(4\gamma)^{1/3}}e^{i(\gamma NS^2-\gamma ^2S^5/10)}\tilde{A}_{32}(S,N-\gamma S^3/3),
-2e^{i(\gamma NS^2-\gamma ^2S^5/10)}\tilde{A}_{32}(S,N-\gamma S^3/3), 
\end{align}
where, adapting the notation of \cite{hewocksmy19}, 
\begin{align}
\tilde{A}_{32}(X,Y) = \int_{\Gamma_{32}}t e^{i(-Yt^2 - Xt^4/2 + 4t^5/(15\gamma^{1/2}))} dt,
\end{align}
with $\Gamma_{32}$ being any contour starting at $t=e^{i9\pi/10}\infty$ and ending at $t=i\infty$. The tilde on $\tilde{A}_{32}$ is included to indicate that the integral $A_{32}$ of \cite{hewocksmy19} has been modified to include a non-trivial amplitude function $F(t)=t$. We note that $\tilde{A}_{32}(X,Y)$ solves the parabolic wave equation \eqref{eq: pwe}. 

In Figure \ref{fig:1} we show a plot of $|A|$ as a function of the curvilinear coordinates $S$ and $N$, and in Figure \ref{fig:2} we show the corresponding approximate Helmholtz equation solution $e^{iks(x,y)}A(S(x,y),N(x,y))$ as a function of the Cartesian coordinates $(x,y)$. 
In Figure \ref{fig:3} we plot $|A|$ as a function of $N$ on the line $S=-10$, alongside the Airy function approximation \eqref{eq: Aincom}.
In Figure \ref{fig:4} we plot $|A|$ as a function of $\hat{K}$ on the lines $S=5$ and $S=10$, accompanied by the approximations \eqref{eq:modAKplus} (for $\hat{K}> 0$, $|\hat{K}|\gg S^{-1/2}$), \eqref{eq:modAKminus} (for $\hat{K}< 0$, $|\hat{K}|\gg S^{-1/2}$), and \eqref{eq:constant} (for $\hat{K}=O(S^{-1/2})$).  This reveals how, as $S$ increases, the wavefield evolves from a beam-like structure to become more like a shadow boundary.   
In all comparisons we see excellent agreement between the asymptotics and numerics.

\section{Conclusions}
\label{conclusions}

Motivated by the Popov problem of finding the outgoing wavefield generated by a whispering gallery wave as it approaches an inflection point, this paper addresses the wavefield that emerges when an incoming Airy function wave is centred on a curve whose curvature decreases to zero.  
An exact solution can be found in terms of the complicated integral (\ref{eq: solforpopov}) whose integrand contains a branch point.  An asymptotic calculation using the stationary phase method reveals that the radiated field is strongest near the tangent at the inflection point but, unlike a Gaussian beam, it has an asymmetric structure whose  largest amplitude is attained near this tangent and eventually resembles a shadow boundary.
We have also computed the integral using the method of \cite{gibhewhuy23} after making a change of variable so as to obtain a representation without branch points.  
This has yielded wave profiles that compare favourably with the asymptotic predictions.

\subsection*{Acknowledgements}
AG and DH acknowledge support from the EPSRC grant EP/V053868/1, and thank the Isaac Newton Institute for Mathematical Sciences, Cambridge, for support and hospitality during the programme \textit{Mathematical theory and applications of multiple wave scattering}, where work on this paper was undertaken. This work was supported by EPSRC grant no EP/R014604/1. 

The authors are grateful to Valery Smyshlyaev and Michael Berry for helpful discussions.

%%%% Bibtex bibliography, with references stored in the refs.bib file
\bibliographystyle{siam}
\bibliography{refs}
\end{document}